\documentclass[draft|showlink]{siamart1116}

\usepackage{amsfonts}
\usepackage{graphicx}
\usepackage{epstopdf}
\usepackage{algorithmic}
\ifpdf
  \DeclareGraphicsExtensions{.eps,.pdf,.png,.jpg}
\else
  \DeclareGraphicsExtensions{.eps}
\fi

\newtheorem{ass}{Assumption}
\newtheorem{rem}{Remark}

\newcommand{\TheTitleABB}{Adapted  $\th$-Scheme and Its Error Estimates for BSDEs} 
\newcommand{\TheTitle}{Adapted  $\th$-Scheme and Its Error Estimates for Backward Stochastic Differential Equations} 
\newcommand{\TheAuthors}{Chol-Kyu Pak, Mun-Chol Kim and Chang-Ho Rim }

\def\R{\mathbb R} \def\P{\mathbb P} \def\F{\mathcal F} \def\d{\partial} \def\th{\theta}
\def\to{\rightarrow}
\def\beq{\begin{equation}} \def\enq{\end{equation}}
\def\beseq{\begin{subequations}}  \def\enseq{\end{subequations}}
\def\beqa{\begin{eqnarray}} \def\enqa{\end{eqnarray}}
\def\non{\nonumber}
\def\BeDef{\begin{definition}} \def\EnDef{\end{definition}}
\def\BeThe{\begin{theorem}} \def\EnThe{\end{theorem}}
\def\BeLem{\begin{lem}} \def\EnLem{\end{lem}}

\numberwithin{theorem}{section}
\numberwithin{equation}{section}
\numberwithin{table}{section}
\numberwithin{figure}{section}

\headers{\TheTitleABB}{\TheAuthors}

\title{{\TheTitle}}

\author{
  Chol-Kyu Pak\thanks{Faculty of Mathematics, Kim Il Sung University, Pyongyang,  Democratic People's Republic of Korea (\email{pck2016217@gmail.com}).}
\and
  Mun-Chol Kim\footnotemark[1]
\and
  Chang-Ho Rim\footnotemark[1]
}

\usepackage{amsopn}

\begin{document}

\maketitle

\begin{abstract}
In this paper we propose a new kind of high order numerical scheme for backward stochastic differential equations(BSDEs). Unlike the traditional $\th$-scheme, we reduce truncation errors by taking $\th$ carefully for every subinterval according to the characteristics of integrands. We give error estimates of this nonlinear scheme and verify the order of scheme through a typical numerical experiment.
\end{abstract}



\section{Introduction}\label{sec:intro}
Let $(\Omega,\F, \P)$ be a probability space, $T>0$ a finite time and $\{\F_t\}_{0\le t \le T}$ a filtration satisfying the usual conditions. Let $(\Omega,\F, \P,\{\F_t\}_{0\le t\le T})$ be a complete, filtered probability space on which a standard $d$-dimensional Brownian motion $W_t=(W_t^1,W_t^2, \cdots ,W_t^d )^T$ is defined and $\F_0$ contains all the $\P$-null sets of $\F$. Let $L^2=L_\F^2 (0,T)$ be the set of all ${\{\F_t\}}$-adapted mean-square-integrable processes.
\par The general form of backward stochastic differential equation (BSDE) is
\beq\label{eq:bsde}
y_t=\xi+\int_t^T{f(s,y_s,z_s )ds}-\int_t^T{z_sdWs}, \quad  t\in[0,T] 
\enq
where the generator $f=f(t,y,z)$ is a vector function valued in $\R^m$  and is $\{\F_t\}$-adapted for each $(y,z)$ and the terminal variable $\xi\in L^2$ is $\F_T$-measurable.
A process $(y_t,z_t ):[0,T]\times\Omega \rightarrow \R^m\times \R^{m\times d}$ is called an $L^2$-solution of the BSDE {\cref{eq:bsde} if it is ${\{\F_t\}}$-adapted, square integrable and satisfies the equation.
\par In 1990, Pardoux and Peng first proved in \cite{Peng90} the existence and uniqueness of the solution of general nonlinear BSDEs and afterwards there has been very active research in this field with many applications.(\cite{Karoui97})
\par In this paper we assume that the terminal condition is a function of $W_T$, i.e. $\xi=\varphi(W_T)$ and the BSDE \cref{eq:bsde} has a unique solution $(y_t,z_t )$. 
It was shown in \cite{Peng91} that the solution $(y_t,z_t)$ of \cref{eq:bsde} can be represented as
\beq\label{eq:eq12}
y_t=u(t,W_t ),    z_t=\nabla_x u(t,W_t ),  \;  \forall t\in[0,T)
\enq
where $u(t,x)$  is the solution of the parabolic partial differential equation
\beq\label{eq:pde}
\frac{\d u}{\d t}+\frac{1}{2} \sum_{i=1}^d \frac{\d^2 u}{\d x_i^2}+f(t,u,\nabla_x u)=0
\enq
with the terminal condition $u(T,x)=\varphi (x)$, and $\nabla_x u$ is the gradient of $u$ with respect to the spatial variable $x$. The smoothness of $u$ depends on $f$ and $\varphi$.
\par Although BSDEs have very important applications in many fields such as mathematical finance and stochastic control, it is well known that it is difficult to obtain analytical solutions except some special cases and there have been many works on numerical methods to get approximate solution. 
A four step algorithm was proposed in \cite{Prot94} to solve a class of more general equations called forward-backward stochastic differential equations (FBSDEs) and in \cite{Chev} a numerical method based on binomial approach was proposed. Besides, there are very interesting numerical methods for BSDEs (\cite{Bou, Gobet, Prot02, Zhang04}).
\par In 2006, Zhao proposed a new kind of numerical method for BSDE in \cite{Zhao06}, it is called $\th$-scheme and was proved to be very effective through many experiments. This $\th$-scheme is simple in form, stable and fairly accurate. 
In \cite{Zhao09} Zhao et al. proved that $\th$-scheme is of first-order when $\th\neq\frac{1}{2}$ and of second-order when $\th=\frac{1}{2}$ (Crank-Nicolson scheme or C-N scheme) in the case where the generator is independent of $z$, and in \cite{Zhao13} they proved the same result in the case of general generators. So we can not expect the high order convergence with this ``traditional"  $\th$-scheme which uses a fixed $\th$ all the time. Afterwards, there have been works on high order scheme and in \cite{Zhao10} a family of multi-step schemes were proposed. But  to one's regret this multi-step scheme is not stable for $z$ and the convergence order was not satisfactory. In \cite{Zhao12} a generalized $\th$-scheme that is more flexible with some additional coefficients was proposed but the convergence order was still about 2 or 3.
\par In this paper we propose a new kind of high order numerical scheme for BSDEs.
This new scheme is similar to the traditional $\th$-scheme and it could achieve high order convergence rate. The main idea is to take $\th$ carefully for every subinterval according to the characteristics of the integrand. We call this scheme ``adapted $\th$-scheme" because we pick various $\th$ at every subinterval according to the integrand. To the best to our knowledge, this kind of scheme has not been proposed before.
\par We consider the case where the generator $f$ is independent of $z$. We assume that $f$ and $\varphi$ are all bounded, smooth enough and their derivatives are also bounded as in \cite{Zhao09}.
\par The rest of this paper is organized as follows. In \Cref{sec:sec2} we explain the idea of the adapted  $\th$-scheme through the approximation of the integral of real functions. In \Cref{sec:sec3}, we propose a new kind of discrete scheme for BSDEs based on the adapted $\th$-scheme. In \Cref{sec:sec4}, we give error estimates of the new scheme theoretically. In \Cref{sec:sec5}, we give a numerical experiment for a typical BSDE to demonstrate the high order convergence of our scheme. In \Cref{sec:conc}, some conclusions are given.
\section{Approximation of integral based on the adapted $\th$-scheme}\label{sec:sec2}
\subsection{The case where the derivatives of integrand are known}\label{sec:sec21}
Let $f(t):[a,b]\to\R$  be  $q$ times continuously differentiable on  $[a,b]$,   $q+1$ times differentiable on  $(a,b)$ and assume   $f,f',\cdots,f^{(q+1)}$are all bounded.
Now we consider the approximation of the integral of  $f(t)$ on $[a,b]$  .
\[I=\int_a^b{f(s)ds}\]
Let  $a=t_0<t_1<\cdots<t_N=b$ be an equidistant partition of  $[a,b]$ and $h=\frac{T}{N}$ .
We approximate the integral on the subinterval  $[t_n,t_{n+1}]$, $I_n=\int_{t_n}^{t_{n+1}}{f(s)ds}$  by
\beq\label{eq:eq21}
\widehat{I}_n(\theta)=[\theta f(t_n)+(1-\theta)f(t_{n+1})]h
\enq
From Taylor expansion we have 
\begin{align}
\Delta I_n(\th)&=\widehat I_n(\th)-I_n=\int_{t_n}^{t_{n+1}}{[\th(f(s)-f(t_n))+(1-\th)(f(s)-f(t_{n+1}))]ds}\non\\
&=\th\int_{t_n}^{t_{n+1}}{[f'(t_n)(s-t_n)+\cdots+f^{(q)}(t_n)\frac{(s-t_n)^q}{q!}]ds}\non\\
&+(1-\th)\int_{t_n}^{t_{n+1}}{[f'(t_{n+1})(s-t_{n+1})+\cdots+f^{(q)}(t_{n+1})\frac{(s-t_{n+1})^q}{q!}]ds}+R_n(\th)\non
\end{align}
\noindent where
\begin{align}
R_n(\th)&=\int_{t_n}^{t_{n+1}}{\left[\th f^{(q+1)}(\alpha_1)\frac{(s-t_n)^{q+1}}{(q+1)!}+(1-\th) f^{(q+1)}(\alpha_2)\frac{(s-t_{n+1})^{q+1}}{(q+1)!}\right]ds}\non\\
&=\th f^{(q+1)}(\alpha_1)\frac{h^{q+2}}{(q+2)!}+(1-\th) f^{(q+2)}(\alpha_2)\frac{(-1)^{q+2}h^{q+2}}{(q+2)!}\non
\end{align}
for some $\alpha_1, \alpha_2\in[t_n,t_{n+1}]$ and we deduce
\begin{align}
\Delta &I_n(\th)-R_n(\th)=\non\\
&=\th\int_{t_n}^{t_{n+1}}{\sum_{k=1}^{q}{f^{(k)}(t_n)\frac{(s-t_n)^k}{k!}}ds}+(1-\th)\int_{t_n}^{t_{n+1}}{\sum_{k=1}^{q}{f^{(k)}(t_{n+1})\frac{(s-t_{n+1})^k}{k!}}ds}\non\\
&=\th\sum_{k=1}^{q}{\frac{[f^{(k)}(t_n)+(-1)^{k+1}f^{(k)}(t_{n+1})]}{(k+1)!}h^{k+1}}-\sum_{k=1}^{q}{\frac{(-1)^{k+1}f^{(k)}(t_{n+1})}{(k+1)!}h^{k+1}}\non.
\end{align}
So if we take $\th$  as
\beq\label{eq:eq22}
\th_n^q=\frac{\sigma_n^q}{\rho_n^q}=\frac{\sum\limits_{k=1}^{q}{\frac{(-1)^{k+1}f^{(k)}(t_{n+1})}{(k+1)!}h^{k+1}}}{\sum\limits_{k=1}^{q}{\frac{[f^{(k)}(t_n)+(-1)^{k+1}f^{(k)}(t_{n+1})]}{(k+1)!}h^{k+1}}}
\enq
then $\Delta I_n(\th_n^q)=R_n(\th_n^q)$. (Note that we assumed that the denominator is not zero.)
The truncation error becomes
\begin{align}
|\Delta I_n(\th_n^q)|=|R_n(\th_n^q)|&\leq \frac{h^{q+2}}{(q+2)!}max\{f^{(q+1)}(\alpha_1),f^{(q+1)}(\alpha_2)\}(2|\th_n^q|+1)\non\\
&\leq(2|\th_n^q|+1)C_{q+1}h^{q+2}\non
\end{align}
where  $C_{q+1}$ is a constant which depends only on the bound of the  $(q+1)$th derivative of  $f$. \\

\BeDef[Validity of subinterval in integral approximation]\label{def:def21}
\\For a constant  $L_\th>0$, the subinterval  $[t_n,t_{n+1}] (0\leq n \leq N-1 )$ is said to be \textbf{valid} if 
\[ \rho_n^q\neq0, \;|\th_n^q|\leq L_\th \] 
where $\rho_n^q,\th_n^q$  are defined in \cref{eq:eq22}\\
\EnDef

Now we take $\th_n=\th_n^q$  for valid subintervals and  $\th_n=\frac{1}{2}$ for invalid ones in the adapted $\th$-scheme. Then the truncation error in valid subintervals satisfies
\beq\label{eq:eq23}
|\Delta I_n|=|\Delta I_n(\th_n^q)|\leq(2L_\th+1)C_{q+1}h^{q+2}
\enq
\noindent  and in the invalid subintervals it is equal to the C-N scheme and satisfies
\[|\Delta I_n|\leq C_2h^3.\] 
So if there are $M$  invalid subintervals, the overall truncation error satisfies
\beq\label{eq:eq24}
|\Delta I|\leq \sum_{n=0}^{N-1}{ |\Delta I_n|}\leq (N-M)(2L_\th+1)C_{q+1}h^{q+2}+MC_2h^3
\enq
We call the above method ``the $q$th order adapted $\th$-scheme" for the approximation of the integral.
\par $\th_n^q$ for $q=1,2,3$ are as follows.\\
- For $q=1$
\beq \label{eq:eq25} \th_n^1=\frac{f'(t_{n+1})}{f'(t_n)+f'(t_{n+1})}\enq
- For $q=2$
\beq \label{eq:eq26} \th_n^2=\frac{3f'(t_{n+1})-f''(t_{n+1})h}{3[f'(t_n)+f'(t_{n+1})]+[f''(t_n)-f''(t_{n+1})]h}\enq
- For $q=3$
\beq \label{eq:eq27} \th_n^3=\frac{12f'(t_{n+1})-4f''(t_{n+1})h+f'''(t_{n+1})h^2}{12[f'(t_n)+f'(t_{n+1})]+4[f''(t_n)-f''(t_{n+1})]h+[f'''(t_{n})+f'''(t_{n+1})]h^2}\enq
\par Now let us discuss under what conditions the subinterval  $[t_n,t_{n+1}]$ is valid.
\\In the case of  $q=1$, if  $f'(t_n)f'(t_{n+1})\geq0$ we clearly have  $\th_n^1\in[0,1]$.
And if  $f'(t_n)f'(t_{n+1})<0$, there exists  $s\in[t_n,t_{n+1}]$ such that  $f'(s)=0$ by intermediate-value theorem, and we deduce that the number of invalid subintervals does not exceed the number of subintervals that have points at which  $f'(t)$ becomes zero.\\
Likewise in the case of  $q>1$,   if
\beq\label{eq:eq28}
\left(\sum_{k=1}^q{\frac{f^{(k)}(t_n)}{(k+1)!}h^{k+1}}\right)\left(\sum_{k=1}^q{\frac{(-1)^{k+1}f^{(k)}(t_{n+1})}{(k+1)!}h^{k+1}}\right)\geq0
\enq
then $\th_n^q\in[0,1]$ and $[t_n,t_{n+1}]$  is valid. From the fact that  $\th_n^q\to \th_n^1$ as $h\to0$ , it would be similar to the case of  $q=1$ when  $h$ is small enough.
\par After all, we could say that the subintervals around zero points of  $f'(t)$ are likely to be invalid. As the validity of subintervals depend on the partition, it is difficult to obtain the general relationship between the numbers of valid ones and invalid ones. But if there are a finite number of zero points of  $f'(t)$ in $[a,b]$ , the ratio of invalid ones to valid ones would be smaller as we increase the size of the partition.
\subsection{The case where the derivatives of integrand are not known}\label{sec:sec22}
As we shall see later, in the case of BSDEs we do not know the precise derivatives of the integrands. So we discuss about the approximation of $\th_n^q$  in \cref{eq:eq22}.
We approximate the derivatives of $f(t)$ by the ones of the Lagrange interpolation polynomial.\\
Assume that $f(t)$ is $q+1$ times differentiable and $f(t_i ),i=0\cdots q$ are given.\\
If we let 
\beq\label{eq:eq29}I_i (t)=\frac{\Pi(t)}{(t-t_i ) \Pi'(t_i )}, \quad  \Pi(t)=\Pi_{i=0}^q(t-t_i) \enq
the Lagrange interpolation polynomial $L(t)$ can be expressed in the form 
\beq\label{eq:eq210}L(t)=\sum_{i=0}^q{I_i(t)f(t_i)} \enq
and the deviation is given by
\beq\label{eq:eq211}L(t)-f(t)=\frac{f^{(q+1)} (\xi)\Pi(t)}{(q+1)!}  ,\quad  t_0\leq\xi\leq t_n.\enq
For $n\leq N-q-1$, we define $L_n (t)$ as the Lagrange interpolation polynomial based on $q+1$ pairs ${(t_{n+k},f(t_{n+k} )):1\leq k\leq q+1}$.

Now we approximate  $f^{(k)}(t_n)$ and  $f^{(k)}(t_{n+1})$ $(n\leq N-q-1)$ in \cref{eq:eq22} as follows. 
\beqa
f^{(k)}(t_n)\approx	\widetilde{f}^{(k)}(t_n)=L_n^{(k)}(t_n) \label{eq:eq212}\\
f^{(k)}(t_{n+1})\approx	\widetilde{f}^{(k)}(t_{n+1})=L_n^{(k)}(t_{n+1})\label{eq:eq213}
\enqa
Then $\widetilde{f}^{(k)}(t_n),\widetilde{f}^{(k)}(t_{n+1})$  can be written as
\beq\label{eq:eq214}  		\widetilde{f}^{(k)}(t_{n})=h^{-k}\sum_{j=1}^{q+1}{t_{kj}^1f(t_{n+j})},
					 \;\widetilde{f}^{(k)}(t_{n+1})=h^{-k}\sum_{j=1}^{q+1}{t_{kj}^2f(t_{n+j})} \enq
where $t_{kj}^1, t_{kj}^2   (1\leq k \leq q,1\leq j \leq q+1)$ are coefficients of  $f(t_{n+j})$ in (\ref{eq:eq212}) and (\ref{eq:eq213}), respectively.

From \cref{eq:eq211}, we get
\beq\label{eq:eq215} \widetilde{f}^{(k)}(t_{n})-f^{(k)}(t_n)=C_1h^{q-k+1},\; \widetilde{f}^{(k)}(t_{n+1})-f^{(k)}(t_n+1)=C_2h^{q-k+1}\enq
where  $C_1,C_2$ are constants that depend only on the bound of $f^{(q+1)}$ .
\par Now we approximate  $\th_n^q$ as follows.

\beq\label{eq:eq216} 
\th_n^q\approx \widetilde{\th}_n^q=\frac{\widetilde{\sigma}_n^q}{\widetilde{\rho}_n^q}=
\frac{{\sum\limits_{k=1}^{q}{\frac{(-1)^{k+1}\widetilde{f}^{(k)}(t_{n+1})}{(k+1)!}h^{k+1}}}}
{\sum\limits_{k=1}^{q}{\frac{[\widetilde{f}^{(k)}(t_{n})+(-1)^{k+1}\widetilde{f}^{(k)}(t_{n+1})]}{(k+1)!}h^{k+1}}}
\enq

If we let $r_j=\sum_{k=1}^{q}{\frac{(-1)^{k+1}t_{kj}^2}{(k+1)!}},  s_j=\sum_{k=1}^q{\frac{t_{kj}^1}{(k+1)!}}$  we have

\beq\label{eq:eq217}
\widetilde{\th}_n^q=\frac{\widetilde{\sigma}_n^q}{\widetilde{\rho}_n^q}=\frac{\sum\limits_{j=1}^{q+1}{r_jf(t_{n+j})}}{\sum\limits_{j=1}^{q+1}{(r_j+s_j)f(t_{n+j})}}\quad   (0\leq n\leq N-q-1)
\enq
and from \cref{eq:eq215} one can see that
\[
|\widetilde{\rho}_n^q-\rho_n^q|\leq C_1h^{q+2},\;|\widetilde{\sigma}_n^q-\sigma_n^q|\leq C_2h^{q+2}
\] 
where  $C_1, C_2$ are constants that depend only on the bound of  $f^{(q+1)}$.
\par So the deviation of $ \widetilde\th_n^q$ from $\th_n^q$  is
\beseq
\begin{align}
|\widetilde\th_n^q-\th_n^q|
&=\left |\frac{\widetilde\sigma_n^q}{\widetilde\rho_n^q}-\frac{\sigma_n^q}{\rho_n^q}\right|=\left |\frac{\rho_n^q(\widetilde\sigma_n^q-\sigma_n^q)+\sigma_n^q(\rho_n^q-\widetilde\rho_n^q)}{\widetilde\rho_n^q\rho_n^q}\right|\non\\
&\leq\left|\frac{\widetilde\sigma_n^q-\sigma_n^q}{\widetilde\rho_n^q}\right|+\left|\th_n^q\frac{\widetilde\rho_n^q-\rho_n^q}{\widetilde\rho_n^q}\right|\leq\left|\frac{1}{\widetilde\rho_n^q}\right|(C_1h^{q+2}+|\th_n^q|C_2h^{q+2})\non \\
&\leq\left|\frac{1}{\widetilde\rho_n^q}\right|(C_1h^{q+2}+|\th_n^q-\widetilde\th_n^q|C_2h^{q+2}+|\widetilde\th_n^q|C_2h^{q+2})\non
\end{align}
\enseq

If we assume that there exist  $L_\th>0$ and $L_\rho>0$  such that $|\widetilde\th_n^q|\leq L_\th$  and $|\widetilde\rho_n^q|^{-1}\leq L_\rho$ for all $n$ , we have 
\[
|\widetilde\th_n^q-\th_n^q|\leq L_\rho(C_1h^{q+2}+|\th_n^q-\widetilde\th_n^q|C_2h^{q+2}+L_\th C_2h^{q+2}) \non
\]
and 
\[
(1-L_\rho C_2 h^{q+2})|\widetilde\th_n^q-\th_n^q|\leq L_\rho(C_1h^{q+2}+L_\th C_2h^{q+2}). \non
\]
Now we can pick $C_0 \in (0,1)$ such that \[ 1-L_\rho C_2 h^{q+2}\geq C_0\] by choosing  $h$ small enough and for this  $C_0$ we deduce 

\addtocounter{equation}{-1} 
\beq\label{eq:eq218}
 |\widetilde\th_n^q-\th_n^q|\leq \frac{L_\rho(C_1+L_\th C_2)}{C_0}h^{q+2}\leq Ch^{q+2}
\enq
where $C$  is a constant that depends only on  $L_\rho, L_\th$ and the bound of  $f^{q+1}$.\\
\BeDef[Validity of subinterval in integral approximation using approximate derivatives]\label{def:def22} 
\\For constants $L_\th>0$ and $ L_\rho>0$ , the subinterval   $[t_n,t_{n+1}](0\leq n \leq N-q-1)$ is said to be \textbf{valid} if 
\beq\label{eq:eq219}
|\widetilde\th_n^q|\leq L_\th, \;|\widetilde\rho_n^q|^{-1}\leq L_\rho
\enq
where  $\widetilde\rho_n^q,\widetilde\th_n^q$ are defined in \cref{eq:eq217}\\
\EnDef 

\addtocounter{equation}{-1} 

\par Now as in \Cref{sec:sec21}, we take $\th_n=\widetilde\th_n^q$  for valid subintervals and  $\th_n=\frac{1}{2}$ for invalid ones.

Then the truncation error in the valid subinterval $[t_n,t_{n+1}]$  becomes
\beseq
\begin{align}
|I_n-\widehat I_n(\widetilde\th^q)|\non
&\leq |I_n-\widehat I_n(\th^q)|+|\widehat I_n(\th^q)-\widehat I_n(\widetilde\th^q)|\non\\
&\leq(2L_\th+1)C_{q+1}h^{q+2}+|f(t)-f(t+h)||\th_q-\widetilde\th_q|h\non\\
&\leq(2L_\th+1)C_{q+1}h^{q+2}+C_1L_\rho(1+L_\th)h^{q+4}\non\\
&\leq Ch^{q+2} \non
\end{align}
\enseq
that is
\beq\label{eq:eq220}
|I_n-\widehat I_n(\widetilde\th^q)|\leq Ch^{q+2}.
\enq
If we assume that the truncation errors for $n>N-q-1$  do not exceed $Ch^{q+2}$  and the number of invalid subintervals is $M$ , the overall error does not exceed
\beq\label{eq:eq221}
(N-M)C_1h^{q+2}+MC_2h^3.
\enq
We give  $\widetilde\th_n^q$ for $q=1,2,3,4$  below. (For simplicity we write  $f(t_{n+j})$ as $f_j$.)
\begin{list}{-}{}
\item For $q=1$
\beq\label{eq:eq222} \widetilde\th_n^1=\frac{f_1-f_2}{2(f_1-f_2)}=\frac{1}{2}\enq
\item For $q=2$
\beq\label{eq:eq223} \widetilde\th_n^2=\frac{11f_1-16f_2+5f_3}{12(2f_1-3f_2+f_3)}\enq
\item For $q=3$
\beq\label{eq:eq224} \widetilde\th_n^3=\frac{31f_1-59f_2+37f_3-9f_4}{24(3f_1-6f_2+4f_3-f_4)}\enq
\item For $q=4$
\beq\label{eq:eq225} \widetilde\th_n^4=\frac{1181f_1-2774f_2+2616f_3-1274f_4+251f_5}{720(4f_1-10f_2+10f_3-5f_4+f_5)}\enq
\end{list}
Note that the case of $q=1$  is equivalent to the C-N scheme.

\subsection{An Example}\label{sec:sec23}
We test the efficiency of the adapted $\th$-scheme through the approximation of the integral of $f(t)=t^3e^{-\left(t-\frac{1}{2}\right)^2}$  on  $[-3, 3]$. We try C-N scheme and the adapted $\th$-scheme of order 2, 3(using \cref{eq:eq223} and \cref{eq:eq224}) for comparison. We set  $L_\rho=1e+8,L_\th=1$ and compare the errors increasing the size of partition from $2^7$ to $2^{12}$. As the size of the partition grows  $\widetilde\rho$ becomes very small, so we take $L_\rho$  large to reduce the number of invalid subintervals. The size of partitions and the number of invalid subintervals in the experiment are shown in \cref{tab:t21}. In the table, we denot by TOT and INV the size of partition and the number of invalid subintervals, respectively.
\begin{table}[tbhp]
\caption{The size of partitions and the number of invalid subintervals}
\label{tab:t21}
\centering
\begin{tabular}{|c|c|c|c|} \hline
No & TOT & INV $(q=2)$ & INV $(q=3)$\\ \hline
1 & 128 & 1 & 0 \\
2 & 256 & 2 & 1 \\
3 & 512 & 1 & 0 \\ 
4 & 1024 & 2 & 1 \\ 
5 & 2048 & 1 & 0 \\ 
6 & 4096 & 2 & 0 \\ 
\hline
\end{tabular}
\end{table}
In this example $f'(t)$  is zero at a point in  $\left\{-1,0,\frac{3}{2}\right\}$ and the invalid subintervals appear around $-1$ and $0$. Note that compared to the size of partitions being increased, the number of invalid subintervals does not exceed 2. (On the other hand we repeat the same experiment increasing the size of partition from $3^4$ to $3^{10}$ and no invalid subintervals appear.)
\cref{fig:fig21} shows the convergence result in \texttt{log} scale. The convergence rate(denoted by CR) is obtained by linear least square fitting. These CR values are consistent with the theoretical result. Especially one can see that the error for 3rd order adapted $\th$-scheme rises a little at the fourth point  and this is because there is an invalid subinterval at Experiment 4.(see \cref{tab:t21}) We perform some more tests for various types of integrands and get similar results.
\begin{figure}[ht]\label{fig:fig21}
\begin{center}
\includegraphics[width=89mm,height=56mm]{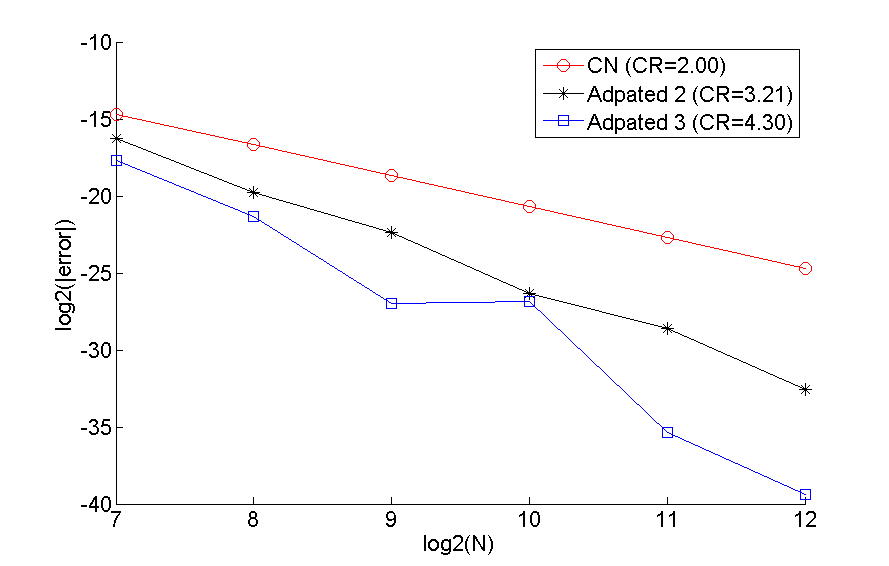}\\
\caption{Errors for approximation of integral using the adapted $\th$-scheme}
\end{center}
\end{figure}

\section{Discrete scheme for BSDEs based on the adapted $\th$-scheme}\label{sec:sec3}
In this section, we propose a discrete scheme for BSDE: 
\beq\label{eq:eq31} y_t=\varphi(W_T)+\int_t^T{f(s,y_s)ds}-\int_t^T{z_sdW_s}, \quad t\in [0,T] \enq
For the sake of simplicity, we assume that everything is one-dimensional, i.e. $m=d=1$ , but all discussions can be generalized to the multi-dimensional case easily.
\par Let $0=t_0<\cdots<t_N=T$  be an equidistant partition of the time interval  $[0,T]$ and $t_{n+1}-t_n=h=\frac{T}{N}$. 
\par In $[t_n, t_{n+1}]$ the BSDE \cref{eq:eq31} can be written as follows.
\beq\label{eq:eq32} y_{t_n}=y_{t_{n+1}}+\int_{t_n}^{t_{n+1}}{f(s,y_s)ds}-\int_{t_n}^{t_{n+1}}{z_sdW_s} \enq
Let $\F_s^{t,x}  (t\leq s\leq T)$ be a $\sigma$-field generated by the Brownian motion $\{x+W_r-W_t,t\leq r\leq s\}$ starting from the time-space point  $(t,x)$ and $E_s^{t,x} [X]:=E[X|\F_s^{t,x}], E_t^x [X]:=E[X|\F_t^{t,x}]$ as in \cite{Zhao13,Zhao09}.
\\Taking $E_{t_n}^x[\cdot]$  to the both sides of \cref{eq:eq32} leads to
\beq\label{eq:eq33} y_{t_n}=E_{t_n}^x[y_{t_{n+1}}]+\int_{t_n}^{t_{n+1}}{E_{t_n}^x[f(s,y_s)]ds}\enq
where the integrand $E_{t_n}^x[f(s,y_s)]$ is a deterministic function of  $s$.
\par Now we introduce the variational equation of \cref{eq:eq31} as follows.
\beq\label{eq:eq34} \nabla y_t=\varphi_x(W_T)+\int_t^T{f_y(s,y_s)\nabla y_sds}-\int_t^T{\nabla z_sdW_s}\enq
where  $\varphi_x,f_y$ are the partial derivatives of $\varphi,f$  with respect to  $x,y$ and $\nabla y_s,\nabla z_s$ are the variations of $y_s,z_s$ with respect to spatial variable  $x$.(See \cite{Zhao09})

From \cref{eq:eq12} we have\[z_t=\nabla y_t\]
and \cref{eq:eq34} can be written as
\beq\label{eq:eq35} z_t=\varphi_x(W_T)+\int_t^T{f_y(s,y_s)z_sds}-\int_t^T{\nabla z_sdW_s} \enq
So we have
\beq\label{eq:eq36} z_{t_n}=z_{t_{n+1}}+\int_{t_n}^{t_{n+1}}{f_y(s,y_s)z_sds}-\int_{t_n}^{t_{n+1}}{\nabla z_sdW_s} \enq
and taking $E_{t_n}^x[\cdot]$  to the both sides leads to
\beq\label{eq:eq37} z_{t_n}=E_{t_n}^x[z_{t_{n+1}}]+\int_{t_n}^{t_{n+1}}{E_{t_n}^x[f_y(s,y_s)z_s]ds}\enq
In \cite{Zhao13,Zhao09} the integrals in \cref{eq:eq33} and \cref{eq:eq37} were replaced by the approximation based on $\th$-scheme to get the discrete scheme, and in \cite{Zhao10} the integrands were replaced by their interpolation polynomials resulting in multi-step scheme.
\par The traditional $\th$-scheme is based on 
\beq\label{eq:eq38} \int_{t_n}^{t_{n+1}}{E_{t_n}^x[f(s,y_s)]ds}=h\left(\th f(t_n,y_{t_n})+(1-\th)E_{t_n}^x[f(t_{n+1},y_{t_{n+1}})]\right)+R_y^n \enq

\begin{multline}\label{eq:eq39} 
\int_{t_n}^{t_{n+1}}{E_{t_n}^x[f_y(s,y_s)z_s]ds}=\\h(\th f_y(t_n,y_{t_n})z_{t_n}+
(1-\th)E_{t_n}^x[f_y(t_{n+1},y_{t_{n+1}})z_{t_{n+1}}])+R_z^n. 
\end{multline}
where $\th$ is a global constant and we have the reference equation as follows. 
\begin{multline}\label{eq:eq310}
\\y_{t_n}=E_{t_n}^x[y_{t_{n+1}}]+h\left(\th f(t_n,y_{t_n})+(1-\th)E_{t_n}^x[f(t_{n+1},y_{t_{n+1}})]\right)+R_y^n\\
z_{t_n}=E_{t_n}^x[z_{t_{n+1}}]+h\left(\th f_y(t_n,y_{t_n})z_{t_n}+(1-\th)E_{t_n}^x[f_y(t_{n+1},y_{t_{n+1}})z_{t_{n+1}}]\right)+R_z^n\\
(0\leq n \leq N-1).
\end{multline}
The discrete scheme based on \cref{eq:eq310} is
\beq\label{eq:eq311}
\left \{ \begin{array}{l}
y^n=E_{t_n}^x[y^{n+1}]+h\left(\th f(t_n,y^n)+(1-\th)E_{t_n}^x[f(t_{n+1},y^{n+1})]\right) \\
z^n=E_{t_n}^x[z^{n+1}]+h\left(\th f_y(t_n,y^n)z^n+(1-\th)E_{t_n}^x[f_y(t_{n+1},y^{n+1})z^{n+1}]\right) \\
y^N=\varphi(W_T) \\
z^N=\frac{d\varphi}{dx}(W_T)
\end{array}\right.
\enq
and it is proved that this scheme achieves the best convergence rate, $2$, when $\th=\frac{1}{2}$ .
\par Here we approximate the integrals in \cref{eq:eq38} and \cref{eq:eq39} using the $q$th order adapted $\th$-scheme which needs  $q+1$ times differentiability of the integrands and the following lemma guarantees this. (\cite{Zhao09}) \\
\begin{lemma}\label{lem:lem31}
Let $\Delta_t^xW_s=x+W_s-W_t$ and $g(s,x),v(s,x),w(s,x)$  be certain functions and $G(s),\widehat G(s)$ be as follows.
\[ G(s)=E_t^x[g(s,v(s,\Delta_t^xW_s))], \widehat G(s)=E_t^x[g(s,v(s,\Delta_t^xW_s))w(s,\Delta_t^xW_s)]\]
If there exists a positive integer  $m$ such that for all $\beta_1, \beta_2$  satisfying $0\leq\beta_1\leq m+1,0\leq\beta_2\leq 2m+1, \beta_1+\beta_2\leq 2m+1$ the derivatives
\[\frac{\partial^{\beta_1+\beta_2}g(s,x)}{\partial^\beta_1s\partial^\beta_2x}, \frac{\partial^{\beta_1+\beta_2}v(s,x)}{\partial^\beta_1s\partial^\beta_2x},
\frac{\partial^{\beta_1+\beta_2}w(s,x)}{\partial^\beta_1s\partial^\beta_2x}
\]
 are continuous and bounded, then $G(s),\widehat G(s)$  are $m$  times continuously differentiable and the derivatives are also bounded.
\end{lemma}
\begin{proof}
Applying Ito's formula to $g(s,v(s,\Delta_t^xW_s))$ and $g(s,v(s,\Delta_t^xW_s))\cdot\\w(s,\Delta_t^xW_s)$  repeatedly, the proof is straightforward. 
\end{proof}
\par So under the assumption that the parameters of \cref{eq:eq31} are smooth enough, the integrands $E_{t_n}^x[f(s,y_s)]$  and  $E_{t_n}^x[f_y(s,y_s)z_s]$ are also smooth enough and we can apply the adapted $\th$-scheme of a proper order.
\par Based on the $q$ th order adapted $\th$-scheme we have the reference equations as follows.
\begin{multline}\label{eq:eq312}
\\y_{t_n}=E_{t_n}^x[y_{t_{n+1}}]+h\left(\widetilde\th_n^y(x) f(t_n,y_{t_n})+(1-\widetilde\th_n^y(x))E_{t_n}^x[f(t_{n+1},y_{t_{n+1}})]\right)+R_y^n\\
z_{t_n}=E_{t_n}^x[z_{t_{n+1}}]+h\left(\widetilde\th_n^z(x) f_y(t_n,y_{t_n})z_{t_n}+(1-\widetilde\th_n^z(x))E_{t_n}^x[f_y(t_{n+1},y_{t_{n+1}})z_{t_{n+1}}]\right)+R_z^n\\
(0\leq n \leq N-1)
\end{multline}
where
\begin{multline}\label{eq:eq313}
\begin{aligned}
R_y^n=\int_{t_n}^{t_{n+1}}E_{t_n}[f(s,y_s)]ds&-h\left(\widetilde\th_n^y(x) f(t_n,y_{t_n})+(1-\widetilde\th_n^y(x))E_{t_n}^x[f(t_{n+1},y_{t_{n+1}})]\right)\\
R_z^n=\int_{t_n}^{t_{n+1}}E_{t_n}[f_y(s,y_s)z_s&]ds-h\big(\widetilde\th_n^z(x) f_y(t_n,y_{t_n})z_{t_n}+\\
&+(1-\widetilde\th_n^z(x))E_{t_n}^x[f_y(t_{n+1},y_{t_{n+1}})z_{t_{n+1}}]\big)
\end{aligned}
\end{multline}
and $\widetilde\th_n^y(x),\widetilde\th_n^z(x)$  are defined as follows. 

\beqa\label{eq:eq314}
\begin{aligned}
\widetilde{\th}_n^y(x)=\frac{\widetilde{\sigma}_n^y(x)}{\widetilde{\rho}_n^y(x)}=\frac{E_{t_n}^x[\sum\limits_{j=1}^{q+1}{r_jf(t_{n+j},y_{t_{n+j}})}]}{E_{t_n}^x[\sum\limits_{j=1}^{q+1}{(r_j+s_j)f(t_{n+j},y_{t_{n+j}})}]}\\
\widetilde{\th}_n^z(x)=\frac{\widetilde{\sigma}_n^z(x)}{\widetilde{\rho}_n^z(x)}=\frac{E_{t_n}^x[\sum\limits_{j=1}^{q+1}{r_jf_y(t_{n+j},y_{t_{n+j}})z_{t_{n+j}}}]}{E_{t_n}^x[\sum\limits_{j=1}^{q+1}{(r_j+s_j)f_y(t_{n+j},y_{t_{n+j}})z_{t_{n+j}}}]}\\
(0\leq n\leq N-q&-1)
\end{aligned}
\enqa
Note that  $\widetilde\th_n^\cdot$ depends on the space point  $x$ because the integrand depends on  $x$, that is, even in the same time interval  $\widetilde\th_n^y$ and $\widetilde\th_n^z$ differ according to the space point. If any of them is not well defined, we use $\frac{1}{2}$  as in \Cref{sec:sec21} and \Cref{sec:sec22}.\\

\BeDef[Validity of subinterval in reference equations]\label{def:def32} 
\\For constants $L_\rho>0$  and $L_\th>0$ , the subinterval  $[t_n,t_{n+1}] (0\leq n \leq N-q-1 )$ is said to be \textbf{valid} if 
\beqa
|\widetilde\rho_n^y(x)|^{-1}\leq L_\rho, |\widetilde\th_n^y(x)|\leq L_\th\non\\
|\widetilde\rho_n^z(x)|^{-1}\leq L_\rho, |\widetilde\th_n^z(x)|\leq L_\th\non
\enqa
for all  $x$, where $\widetilde\rho_n^y, \widetilde\th_n^y,\widetilde\rho_n^z,\widetilde\th_n^z$ are defined in \textnormal{(\ref{eq:eq314})}.\\
\EnDef
From this reference equation we can get a discrete scheme of (3.1) as follows
\beqa\label{eq:eq315}
\begin{aligned}
y^n=E_{t_n}^x[y^{n+1}]+h(\widehat\th_n^y f(t_n,y^n)+(1-\widehat\th_n^y)E_{t_n}^x[&f(t_{n+1},y^{n+1})]) \\
z^n=E_{t_n}^x[z^{n+1}]+h(\widehat\th_n^z f_y(t_n,y^n)z^n+(1-\widehat\th_n^z)E_{t_n}^x[&f_y(t_{n+1},y^{n+1})z^{n+1}]) \\
&0\leq n\leq N-q-1
\end{aligned}
\enqa
where
\beqa\label{eq:eq316}
\begin{aligned}
\widehat{\th}_n^y(x)=\frac{\widehat{\sigma}_n^y(x)}{\widehat{\rho}_n^y(x)}=\frac{E_{t_n}^x[\sum\limits_{j=1}^{q+1}{r_jf(t_{n+j},y^{n+j})}]}{E_{t_n}^x[\sum\limits_{j=1}^{q+1}{(r_j+s_j)f(t_{n+j},y^{n+j})}]}\\
\widehat{\th}_n^z(x)=\frac{\widehat{\sigma}_n^z(x)}{\widehat{\rho}_n^z(x)}=\frac{E_{t_n}^x[\sum\limits_{j=1}^{q+1}{r_jf_y(t_{n+j},y^{n+j})z^{n+j}}]}{E_{t_n}^x[\sum\limits_{j=1}^{q+1}{(r_j+s_j)f_y(t_{n+j},y^{n+j})z^{n+j}}]}\\
(0\leq n\leq N-q&-1)
\end{aligned}
\enqa
Here we assume that we have approximations $(y^j,z^j)_{N-q\leq j \leq N}$  using any other numerical methods , for example C-N scheme. 
\par We call the discrete scheme (\ref{eq:eq315}) the \textbf{``($q$th order) adapted $\th$-scheme"} for BSDE \cref{eq:eq31}. 
\\
\BeDef[Validity of subinterval in the discrete scheme]\label{def:def33} 
\\For constants $L_\rho>0$  and $L_\th>0$ , the subinterval  $[t_n,t_{n+1}] (0\leq n \leq N-q-1 )$ is said to be \textbf{valid} if 
\beqa
|\widehat\rho_n^y(x)|^{-1}\leq L_\rho, |\widehat\th_n^y(x)|\leq L_\th\non\\
|\widehat\rho_n^z(x)|^{-1}\leq L_\rho, |\widehat\th_n^z(x)|\leq L_\th\non
\enqa
for all  $x$, where $\widehat\rho_n^y, \widehat\th_n^y,\widehat\rho_n^z,\widehat\th_n^z$ are defined in \textnormal{(\ref{eq:eq316})}.\\
\EnDef

We note that in invalid subintervals, we use $\frac{1}{2}$ for $\th$.
\\In the case of $q=2$ , from \cref{eq:eq223},  $\widehat\th_n^y(x),\widehat\th_n^z(x)$  can be written as follows.
\beqa\begin{aligned}
&\widehat\th_n^y(x)=\frac{E_{t_n}^x[11f(t_{n+1},y^{n+1})-16f(t_{n+2},y^{n+2})+5f(t_{n+3},y^{n+3})]}{12E_{t_n}^x[2f(t_{n+1},y^{n+1})-3f(t_{n+2},y^{n+2})+f(t_{n+3},y^{n+3})]}\non\\
&\widehat\th_n^z(x)=\frac{E_{t_n}^x[11f_y(t_{n+1},y^{n+1})z^{n+1}-16f_y(t_{n+2},y^{n+2})z^{n+2}+5f_y(t_{n+3},y^{n+3})z^{n+3}]}{12E_{t_n}^x[2f_y(t_{n+1},y^{n+1})z^{n+1}-3f_y(t_{n+2},y^{n+2})z^{n+2}+f_y(t_{n+3},y^{n+3})z^{n+3}]}\non
\end{aligned}\enqa

\begin{rem}\label{rem:rem31}
The adapted $\th$-scheme is similar to the multistep scheme proposed in \cite{Zhao10} that uses approximations at several points, but our new scheme is nonlinear and is stable for both $y$  and $z$  assuming that every subinterval is valid.
\par In fact, the scheme (\ref{eq:eq315}) can be written as 
\[
\left\{\begin{array}{l}
y^n=E_{t_n}^x[y^{n+1}]+hF(t_n,y^n,y^{n+1},\cdots,y^{n+q+1})\\
z^n=E_{t_n}^x[z^{n+1}]+hG(t_n,y^n,y^{n+1},\cdots,y^{n+q+1},z^n,z^{n+1},\cdots,z^{n+q+1})
\end{array}
\right.
\]
and one can see that under the assumption that every subinterval is valid  $F$ and $G$ are Lipschitz continuous.
Furthermore, if $f\equiv0$, $F=0$ and $G=0$. So this scheme is stable from the theory of numercal ODEs.
\end{rem}

\section{Error estimates of adapted $\th$-scheme}\label{sec:sec4}
In this section we give error estimates of the adapted $\th$-scheme proposed in \Cref{sec:sec3}.
First we make some assumptions as follows.\\
\begin{ass}\label{as1}
The functions $\varphi$ and $f$ in \cref{eq:eq31} are bounded, smooth enough with bounded derivatives.
\end{ass}
\begin{ass}\label{as2}
For certain constants $L_\rho>0$  and $L_\th>0$ , every subinterval $[t_n,t_{n+1}]$  is valid in the sense of \Cref{def:def32} and \Cref{def:def33}.\\
\end{ass}
We need the \Cref{as2} only for simplicity and in the case where there are   invalid intervals we could get results similar to \cref{eq:eq221}. \\

\begin{lemma}\label{lem:lem41}
Let $R_y^n,R_z^n$  be truncation errors defined in \textnormal{(\ref{eq:eq313})}. Under the \Cref{as1} and \textnormal{\ref{as2}} we have
\[|R_y^n|\leq Ch^{q+2}, |R_z^n|\leq Ch^{q+2}\] 
where $C$  is a constant that depends only on $L_\rho, L_\th, T$ and bounds of $f,\varphi$ and their derivatives.
\end{lemma}
\begin{proof}
 It can be easily proved using Taylor expansion as in \Cref{sec:sec2}.
\end{proof}
\bigskip
\par Let  $y_t,z_t(t\in[0,T])$ and  $y^n,z^n(0\leq n \leq N-q-1 )$ be solutions of BSDE \cref{eq:eq31} and the discrete scheme (\ref{eq:eq315}) respectively.
Let $e_y^n=y_{t_n}-y^n, e_z^n=z_{t_n}-z^n$ and $e_{\th_y}^n=\widetilde\th_n^y-\widehat\th_n^y ,e_{\th_z}^n=\widetilde\th_n^z-\widehat\th_n^z$ for  $n=0\cdots N-q-1$.\\

\begin{lemma}\label{lem:lem42}
Under the \Cref{as1} and \textnormal{\ref{as2}} the following estimate holds true:
 \[
|e_{\th_y}^n|\leq C\sum_{j=1}^{q+1}{E_{t_n}^x[|e_y^{n+j}|]}
\]
where $C$  is a constant that depends only on  $L_\rho, L_\th, T$ and bounds of $f$ and their derivatives.
\end{lemma}
\begin{proof}
From (\ref{eq:eq314}) and (\ref{eq:eq316}) we have
\begin{eqnarray}\begin{aligned}
\widetilde\rho_y^n&=E_{t_n}^x[\sum_{j=1}^{q+1}{(r_j+s_j)f(t_{n+j},y_{t_{n+j}})}], \widetilde\sigma_y^n=E_{t_n}^x[\sum_{j=1}^{q+1}{r_jf(t_{n+j},y_{t_{n+j}})}]\non\\
\widehat\rho_y^n&=E_{t_n}^x[\sum_{j=1}^{q+1}{(r_j+s_j)f(t_{n+j},y^{n+j})}], \widehat\sigma_y^n=E_{t_n}^x[\sum_{j=1}^{q+1}{r_jf(t_{n+j},y^{n+j})}]\non
\end{aligned}\end{eqnarray}
and
\begin{eqnarray}\begin{aligned}
|\widetilde\rho_y^n-\widehat\rho_y^n|&\leq L_fE_{t_n}^x[\sum_{j=1}^{q+1}{|r_j+s_j||e_y^{n+j}|}]\non\\
|\widetilde\sigma_y^n-\widehat\sigma_y^n|&\leq L_fE_{t_n}^x[\sum_{j=1}^{q+1}{|r_j||e_y^{n+j}|}]\non
\end{aligned}\end{eqnarray}
where $L_f$ is a Lipschitz constant of $f$.
\par So if we set $\gamma=max_{j=1\cdots q+1}\{|r_j|+|s_j|\}$  we have
\beqa
|\widetilde\rho_y^n-\widehat\rho_y^n|\leq \gamma L_f E_{t_n}^x[\sum_{j=1}^{q+1}{|e_y^{n+j}|}]\non\\
|\widetilde\sigma_y^n-\widehat\sigma_y^n|\leq \gamma L_f E_{t_n}^x[\sum_{j=1}^{q+1}{|e_y^{n+j}|}]\non
\enqa
and using the assumptions we can deduce that
\begin{eqnarray}\begin{aligned}
|e_{\th_y}^n|=|\widetilde\th_y^n-\widehat\th_y^n|&\leq \left|\widetilde\th_y^n-\frac{\widehat\sigma_y^n}{\widetilde\rho_y^n}\right|+\left|\widehat\th_y^n-\frac{\widehat\sigma_y^n}{\widetilde\rho_y^n}\right|\non\\
&\leq |\widetilde\rho_y^n|^{-1}\left(|\widetilde\sigma_y^n-\widehat\sigma_y^n|+|\widehat\th_y^n||\widetilde\rho_y^n-\widehat\rho_y^n|\right)\non\\
&\leq L_\rho\left(|\widetilde\sigma_y^n-\widehat\sigma_y^n|+L_\th |\widetilde\rho_y^n-\widehat\rho_y^n|\right)\non\\
&\leq \gamma L_\rho(1+L_\th)L_f E_{t_n}^x[\sum_{j=1}^{q+1}{|e_y^{n+j}|}]\leq CE_{t_n}^x[\sum_{j=1}^{q+1}{|e_y^{n+j}|}]\non
\end{aligned}\end{eqnarray}
which completes the proof.
\end{proof}
\bigskip
\BeThe\label{theo1}
Suppose that  \Cref{as1} and \textnormal{\ref{as2}} hold and that the initial approximation satisfies
 \[\max\limits_{N-q\leq n \leq N}E[|y_{t_n}-y^n|]=O(h^{q+1})\]
Then for sufficiently small time steps  $h$, it holds that
\[\sup\limits_{0\leq n \leq N}E[|y_{t_n}-y^n|]\leq Ch^{q+1}\]
where $C$  is a constant that depends only on $L_\rho, L_\th, T$ and bounds of $f,\varphi$ and their derivatives.
\EnThe
\begin{proof}
For $0\leq n \leq N-q-1$, from (\ref{eq:eq312}) and (\ref{eq:eq315}) we have
\beqa\begin{aligned}
e_y^n&=E_{t_n}^x[e_y^{n+1}]+h\widehat\th_n^y(f(t_n,y_{t_n})-f(t_n,y^n))+hf(t_n,y_{t_n})(\widetilde\th_n^y-\widehat\th_n^y)\non\\
&+h(1-\widehat\th_n^y)E_{t_n}^x[f(t_{n+1},y_{t_{n+1}})-f(t_{n+1},y^{n+1})]+hE_{t_n}^x[f(t_{n+1},y_{t_{n+1}})](\widehat\th_n^y-\widetilde\th_n^y)+R_y^n\non
\end{aligned}\enqa
Let $L_f$ be the Lipschitz constant of $f$, then from the assumptions and \Cref{lem:lem42} we deduce 
\beqa\begin{aligned}
|e_y^n|&\leq E_{t_n}^x[|e_y^{n+1}|]+hL_\th L_f|e_y^n|+h(1+L_\th)L_fE_{t_n}^x[|e_y^{n+1}|]+2hC_0|e_{\th_y}^n|+|R_y^n|\non\\
&\leq E_{t_n}^x[|e_y^{n+1}|]+hL_\th L_f|e_y^n|+h(1+L_\th)L_fE_{t_n}^x[|e_y^{n+1}|]+2hC_0C_1\sum_{j=1}^{q+1}{E_{t_n}^x[|e_y^{n+j}|]}+|R_y^n|\non\\
&\leq E_{t_n}^x[|e_y^{n+1}|]+hC_2\sum_{j=n}^{n+q+1}{E_{t_n}^x[|e_y^{j}|]}+C_3h^{q+2}
\end{aligned}\enqa
where $C_0$   is a constant that depends only on the bound of  $f$, $C_1$  is a constant determined from \Cref{lem:lem42}, $C_2=(1+L_\th)L_f+2C_0C_1$  and $C_3$  is a constant determined from \Cref{lem:lem41}.\\
Likewise, we have
\beqa\begin{aligned}
|e_y^n|\leq E_{t_n}^x[|e_y^{n+1}|]+hC_2&\sum_{j=n}^{n+q+1}{E_{t_n}^x[|e_y^{j}|]}+C_3h^{q+2}\non\\
|e_y^{n+1}|\leq E_{t_n}^x[|e_y^{n+2}|]+hC_2&\sum_{j=n+1}^{n+q+2}{E_{t_n}^x[|e_y^{j}|]}+C_3h^{q+2}\non\\
&\vdots\non\\
|e_y^{N-q-1}|\leq E_{t_n}^x[|e_y^{N-q}|]+hC_2&\sum_{j=N-q-1}^{N}{E_{t_n}^x[|e_y^{j}|]}+C_3h^{q+2}\non
\end{aligned}\enqa
Adding up the above inequalities gives
\beqa\begin{aligned}
|e_y^n|&\leq E_{t_n}^x[|e_y^{N-q}|]+hC_2q\sum_{j=n}^{N}{E_{t_n}^x[|e_y^{j}|]}+NC_3h^{q+2}\non\\
&\leq hC_4\sum_{j=n}^{N}{E_{t_n}^x[|e_y^{j}|]}+C_5h^{q+1}
\end{aligned}\enqa
where we use the assumption on the initial values and  $N=\frac{T}{h}$.
So we have
 \[
|e_y^n|\leq \frac{hC_4}{1-hC_4}\sum_{j=n+1}^{N}{E_{t_n}^x[|e_y^{j}|]}+\frac{C_5}{1-hC_4}h^{q+1}
\]
Now for sufficiently small time steps $h$  that satisfy $1-hC_4>C_6$  for a fixed constant  $C_6\in(0,1)$, we have 
 \[
|e_y^n|\leq hC\sum_{j=n+1}^{N}{E_{t_n}^x[|e_y^{j}|]}+Ch^{q+1}
\]
where $C$ is a generic constant that does not depend on the partition.
\par Let $\zeta_n=hC\sum_{j=n}^N{E_{t_n}^x[|e_y^j|]}+Ch^{q+1}$ then $|e_y^n|\leq \zeta_{n+1}$ and we deduce
\beqa\begin{aligned}
\zeta_n=hC&\sum_{j=n}^N{E_{t_n}^x[|e_y^j|]}+Ch^{q+1}=hC|e_y^n|+\zeta_{n+1}\leq(1+hC)\zeta_{n+1}\non\\
&\leq\cdots\leq(1+hC)^{N-q-n}\zeta_{N-q}\leq\left(1+\frac{TC}{N}\right)^N\zeta_{N-q}\leq e^{TC}\zeta_{N-q}\non
\end{aligned}\enqa
From the assumptions on the initial values we have
\[
E\zeta_{N-q}=hC\sum_{j=N-q}^N{E[|e_y^j|]}+Ch^{q+1}\leq Ch^{q+1}
\]
which says
\[
E|e_y^n|\leq Ch^{q+1}
\]
and the proof is completed.\\
\end{proof}
To get error estimates for $z$ we introduce the following lemma (see \cite{Prot02} for details) to show the boundedness of $z_t$.\\
\begin{lemma}\label{lem:lem44}
Let $z_t$  be the solution of \cref{eq:eq31}. Then under \Cref{as1}, the following estimate holds true: 
 \[\sup\limits_{t\in[0,T]}E[z_t]\leq C\]
where $C$  is a constant depending only on $T$ , upper bounds of the functions $f$ and $\varphi$, and their derivatives.\\
\end{lemma}
Using \Cref{theo1} and \Cref{lem:lem44}, we repeat the procedure of the proof of \Cref{lem:lem42} to get the following estimates.\\
\begin{lemma}\label{lem:lem45}
Under \Cref{as1} and \ref{as2}, the following estimate holds true:
\[
|e_{\th_z}^n|\leq C \sum_{j=1}^{q+1}{E_{t_n}^x[|e_z^{n+j}|]}+Ch^{q+1}
\]
where $C$  is a constant depending only on $L_\rho, L_\th$ , upper bounds of the function $f$ and  its derivatives.\\
\end{lemma}
Using \Cref{lem:lem45}, we repeat the procedure of proof for \Cref{theo1} to get the following estimates.\\

\BeThe\label{theo2} Suppose that \Cref{as1} and \ref{as2} hold and that the initial values  $y^n,z^n (N-q\leq n \leq N )$ satisfy
\[
\max\limits_{N-q\leq n \leq N}E[|y_{t_n}-y^n|]=O(h^{q+1}),\max\limits_{N-q\leq n \leq N}E[|z_{t_n}-z^n|]=O(h^{q+1})
\] 
 Then for sufficiently small time steps  $h$, it holds that
\[
\sup\limits_{0\leq n \leq N}E[|z_{t_n}-z^n|]\leq Ch^{q+1}.
\]
where $C$  is a constant that depends only on $L_\rho, L_\th, T$  and bounds of $f,\varphi$ and their derivatives.
\EnThe
\section{A Numerical Experiment}\label{sec:sec5}
In this section, we present a typical numerical experiment to demonstrate the effect of the adapted $\th$-scheme. We approximate the conditional expectation by Gauss-Hermite quadrature which has been proved to be efficient in many works, and get estimates at nongrid space points by Lagrange interpolation. We take uniform partitions in both time and space in the experiment and choose the space step size $\Delta x$  so that the local error in space to be balanced with the local error in time. When the  $q$th order adapted $\th$-scheme is used, the local error in time is $O(h^{q+2})$  and the local error in space from $r$th order polynomial interpolation is  $O\left((\Delta x)^{r+1}\right)$. So we set $\Delta x=h^{\frac{q+2}{r+1}}$.  We set the number of the Gauss-Hermite quadrature points to be big enough, 8, so that the error contributed by the use of the Gauss-Hermite quadrature rule does not affect measurements of the convergence rate CR.
\par The BSDE to approximate is as follows. (See \cite{Zhao09}.)
\beq\label{eq:eq51}
\left\{\begin{array}{l}
-dy_t=(-y_t^3+2.5y_t^2-1.5y_t)dt-z_tdW_t\\
y_T=\frac{exp(W_T+T)}{exp(W_T+T)+1}
\end{array}\right.
\enq
The analytic solution of this BSDE is given as
\beq\label{eq:eq52}
\left\{\begin{array}{l}
y_t=\frac{exp(W_t+t)}{exp(W_t+t)+1}\\
z_t=\frac{exp(W_t+t)}{\left(exp(W_t+t)+1\right)^2}
\end{array}\right.
\enq
and the exact solution at $t=0$  is  $\left(\frac{1}{2},\frac{1}{4}\right)$.
\par We set  $T=1$ and measure the error at $t=0$  increasing the size of partitions from $2^3$ to $2^7$. For comparison we try C-N scheme and the adapted $\th$-scheme of orders 2, 3, 4 and the errors for $y$ and $z$ are shown in \Cref{tbl51} and \Cref{tbl52}, respectively.(In the table we denote the $q$th order adapted $\th$-scheme by 'Ada $q$'.) The convergence rate CR is obtained by using linear least squares fitting. The convergence rate of the $q$ th order adapted $\th$-scheme is about  $q+1$ theoretically and the experiment result is consistent with the theoretical ones. Note that both the errors for $y$ and $z$ converge at almost the same rate. We set $L_\th=10,L_\rho=1e+30$  in the experiment,  we set $L_\rho$  large enough to reduce the number of invalid subintervals.
\begin{table}[tbhp]\label{tbl51}
\caption{Errors and convergence rates for $y$}
\centering
\begin{tabular}{|c|c|c|c|c|c|c|} \hline
& $N=8$ & $N=16$ & $N=32$ & $N=64$ & $N=128$ & CR\\ \hline
C-N &8.077e-05&2.041e-05&5.146e-06&1.304e-06&3.323e-07& 1.981\\ \hline
Ada 2&6.086e-06&8.907e-07&1.311e-07&1.693e-08&2.210e-09&2.857\\ \hline
Ada 3&3.010e-07&3.327e-08&3.877e-09&2.254e-10&1.985e-11&3.498\\ \hline
Ada 4&2.609e-07&2.108e-09&3.476e-09&2.311e-10&4.450e-13&4.151\\ \hline
\end{tabular}
\end{table}
\begin{table}[tbhp]\label{tbl52}
\caption{Errors and convergence rates for $z$}
\centering
\begin{tabular}{|c|c|c|c|c|c|c|} \hline
& $N=8$ & $N=16$ & $N=32$ & $N=64$ & $N=128$ & CR\\ \hline
C-N&1.124e-04&2.793e-05&6.968e-06&1.723e-06&4.243e-07&2.011\\ \hline
Ada 2&3.232e-05&5.536e-06&6.651e-07&1.009e-07&1.298e-08&2.834\\ \hline
Ada 3&1.226e-05&1.498e-06&8.835e-08&6.215e-09&4.275e-10&3.753\\ \hline
Ada 4&4.516e-06&1.409e-07&1.363e-08&3.926e-10&1.841e-11&4.429\\ \hline
\end{tabular}
\end{table}

\section{Conclusions}\label{sec:conc}
In this paper we proposed a new kind of high-order numerical scheme for BSDEs, called the adapted $\th$-scheme. Unlike the well-known traditional $\th$-scheme, it reduces truncation errors by taking $\th$ adaptively for every subinterval according to the characteristics of the integrand. We gave error estimates of this scheme in the case where the generator $f$  is independent of  $z$ and verified the efficiency of our scheme through a typical numerical experiment. This new scheme is similar to the multistep scheme proposed in \cite{Zhao10} but this scheme is nonlinear and the stability is guaranteed under some reasonable assumptions. The main con of the adapted $\theta$ scheme is that it is not clear for which types of functions the time intervals would be valid. But we think the idea of the adapted $\theta$ scheme would still be useful for other numerical fields. 
\bibliographystyle{siamplain}

\end{document}